\documentclass[12pt]{article}
\usepackage{amsfonts}
\usepackage{amsbsy}
\usepackage{graphicx}
\usepackage{latexsym}
\newcommand{\Th}{{\Theta}}

\newcommand{\dvo}{d{\vec \omega}}

\newcommand{\al}{\alpha}

\newcommand{\de}{\delta}

\newcommand{\pl}{\partial}
\newcommand{\pr}{\prime}

\newcommand{\vo}{{\vec {\omega}}}
\newcommand{\vd}{{\vec {D}}}
\newcommand{\vb}{{\vec {L}}}

\newcommand{\s}{\sigma}
\newcommand{\sm}{\sigma_{-}}

\newcommand{\be}{\begin{equation}}
\newcommand{\ee}{\end{equation}}
\newcommand{\bq}{\begin{eqnarray}}
\newcommand{\eq}{\end{eqnarray}}

\newcommand{\ba}{\begin{array}}
\newcommand{\ea}{\end{array}}

\newcommand{\Rg}{\mathfrak{R}}
\newcommand{\Ra}{R^{(\alpha)}}
\newcommand{\Rb}{R^{(\beta)}}
\newcommand{\ra}{r^{(\alpha)}}
\newcommand{\rb}{r^{(\beta)}}

\newcommand{\cpo}{{\mathbb {CP}}^1}
\newcommand{\Co}{{\mathbb C}}
\newcommand{\cz}{{\mathbb Z}}

\newcommand{\pdk}{\partial_{\de_k}}
\newcommand{\pdj}{\partial_{\de_j}}

\newcommand{\bt}{\beta}
\newcommand{\ca}{{\cal A}}
\newcommand{\om}{\omega}

\newcommand{\pz}{\partial_z}

\renewcommand{\wp}{w_{+}}

\newcommand{\st}{{\sigma_3}}
\newcommand{\stp}{{\sqrt {2\pi i}}}
\newcommand{\ty}{{\tilde Y}}

\newcommand{\hy}{{\hat Y}}

\newcommand{\vP}{{\vec {\Psi}}}
\renewcommand{\a}{{\mathsf{a}}}
\renewcommand{\b}{{\mathsf{b}}}

\newcommand{\tr}{{\rm tr}}
\newcommand{\Om}{{\Omega}}
\begin{document}

\title{\bf A Riemann-Hilbert Approach to the
Akhiezer Polynomials}
\author{Yang Chen$^{{\dag,\dag\dag}}$\\
        Department of Mathematics\\
        University of Wisconsin-Madison\\
        480 Lincoln Drive\\
        Madison, WI 53706, USA\\
        Alexander R. Its$^{*}$\\
        Department of Mathematical Sciences\\
        Indiana University Purdue University Indiana\\
        402 North Blackford Street\\
        Indianapolis, IN 46202-3216,USA }
 \date{01/09/2003}

\maketitle

\begin{abstract}
In this paper, we study those polynomials, orthogonal with respect to a
particular weight, over the unioin of disjoint intervals,
first introduced by N. I. Akhiezer, via a
reformulation as a matrix factorization or Riemann-Hilbert
problem. This approach complements the method proposed in a
previous paper, that involves the construction of
a certain meromorphic function on a
hyperelliptic Riemann surface. The method described
here is based on the general Riemann-Hilbert scheme of
the theory of integrable systems and will
enable us to derive, in a very strightforward way,
the relevant system of Fuchsian differential equations for the
polynomials and the associated system of the Schlesinger
deformation equations for certain quantaties involving the
corresponding recurrence coefficients. Both of these equations
were obtained earlier by A. Magnus. In our approach,
however, we are able to go beyond Magnus's results by actually
solving the equations in terms of the Riemann $\Theta$-functions.  We also
show that the related Hankel determinant can be interpreted as the
relevant ${\cal {\tau}}-$ function.

\end{abstract}
\noindent
$^{\dag}$ ychen@ic.ac.uk,
$^{\dag\dag}$ Address as of 01/01/03:
Department of Mathematics,
Imperial College, 180 Queen's Gates, London SW7 2BZ, UK.\\
\noindent
$^{*}$ itsa@math.iupui.edu\\
\vfill\eject
\noindent
{\bf Acknowledgment}

The first author should like to thank the
of Department of Mathematics, University of Wisconsin-Madison, for the
kind hospitality in hosting him and the EPSRC for a Oversea Travel Grant
that made this endeavour possible. The second author was supported in part
by NSF Grant DMS-0099812 and by Imperial College of the University of London
via a EPSRC Grant. The final part of this project was done when
he was visiting Institut de  Math\'ematique de l' Universit\'e de Bourgogne,
and the support during his stay there is gratefully acknowledged.
\vskip .4cm
\setcounter{equation}{0}

\section{Introduction}

 The Chebyshev polynomials are those monic polynomials characterised by the
property that ${\rm max}|\pi_n(x)|,\;\;x\in[-1,1],\;$ is
as small as possible. Indeed, it is also known that $\pi_n$
is orthogonal with respect to $\frac{1}{{\pi\sqrt {1-x^2}}}$
over $[-1,1].$ The polynomials $\pi_n$---the Chebyshev
polynomials of the first kind---which satisfy a constant
coefficients three term recurrence
relations, can be thought of as the ``Hydrogen Atom'' model
of those polynomials orthogonal over $[-1,1].$ These play a fundamental
role in the large $n$ asymptotics of the
Bernstein-Szeg\"o polynomials which are
orthogonal with respect a ``deformed'' Chebyshev weight,
$p(x)/{\sqrt {1-x^2}},$ over $[-1,1],$ where $p(x)$ is strictly
positive, absolutely continuous and satisfies the Szeg\"o condition
\cite{Sz}
$$
\int_{-1}^{1}\frac{\ln p(x)}{{\sqrt {1-x^2}}}dx>-\infty.
$$
Many years ago N. I. Akhiezer and also Yu. Ya. Tomchuk \cite{Ak}, \cite{AkT},
\cite{Tom} considered
a generalization of the Chebyshev polynomials, where the interval
of orthogonality is a union of disjoint intervals henceforth
denoted as
\be
E:=(\bt_0,\al_1)\cup(\bt_1,\al_2)
\cup\cdots\cup(\bt_g,\bt_{g+1}).
\ee
For comparison with those of Akhiezer, we assume here $\bt_0=-1,$ and
$\bt_{g+1}=1.$ For later convenience, when the end points
become independent variables we shall adopt the convention,
\bq
(\al_1,\al_2,...,\al_g,\;\bt_0,\bt_1,...,\bt_{g+1})
\longrightarrow
(\de_1,\de_2,...,\de_{g+1},\de_{g+2},...,\de_{2g+2}).
\eq
 Let
\bq
w(z):=\frac{i}{\pi}
{\sqrt {\frac{\Pi_{j=1}^{g}(z-\al_j)}{\Pi_{j=0}^{g+1}(z-\bt_j)}}}
,\eq
be defined in the $\cpo\setminus E.$ The
multi-interval analog of the Chebyshev weight is
\be
w_{+}(t)=\frac{1}{\pi}
{\sqrt {\frac{\Pi_{j=1}^{g}(t-\al_j)}
{(\bt_{g+1}-t)(t-\bt_0)\Pi_{j=1}^{g}(t-\bt_j)}}}\:>0,\quad t\in E,
\ee
and is obtained from the continuation $w(z)$ to the top of the cut,
$E.$ The generalized Chebyshev or Akhiezer polynomials $P_n$
are monic polynomials orthogonal with respect to $w_+,$ i.e.,
\bq \label{ortog}
\int_{E}P_{m}(x)P_{n}(x)w_+(x)dx=h_n\de_{m,n},
\eq
where $h_n$ is the square of the $L^2$ norm.

 In the construction of the Bernstein-Szeg\"o asymptotics over
$E$, for polynomials orthogonal with respect to the weight
$p(t)\wp(t),$ where $p$ is an absolutely continuous positive function,
exact information on $P_n$ would be required. This would entail
the solution of the ``Hydrogen Atom'' problem in the multiple interval
situation. In the case of two intervals, $[-1,\al]\cup[\bt,1],$ $P_n$
was constructed by Akhiezer with an innovation which we would
now recognise as the Baker-Akhiezer function, associated with the
discrete Schr\"odinger equation, namely, the three term recurrence relations,
where the degree of the polynomials $n$ is the
``coordinates'', and $z$ is spectral variable. Akhiezer,
based his construction on the conformal mapping of a doubly
connected domain, with the aid of the Jacobian elliptic functions, as
a demonstration for his students, the applications of elliptic
functions \cite{Akbook}. It is not at all clear how the conformal
mapping could be adepted to handle the situation when there
are more then two intervals. In the early 1960's, Akhiezer
and also with Tomchuk published several very short and very deep
papers regarding the Bernstein-Szeg\"o asymptotics.
Akhiezer and Tomchuk gave a description of $P_n$ and $Q_n$ (the
second solution of the recurrence relations) with the aid of theory
of Hyperelliptic integrals in terms of a cerian Abelian integral of 
the third kind.
However, certain unknown
points on Riemann surface appear in this
respresentation, later circumvented in \cite{Chen}.

In a recent work of A. P. Magnus \cite{mag}, a general class of
semi-classical orthogonal polynomials, which includes the
Akhiezer polynomials $P_{n}$, was introduced and shown that
these polynomials satisfy a certain system of linear Fuchsian equations.
It was also demonstrated there that the recurrence
coefficients, as functions of the natural parameters of
the semi-classical weights, obey the nonlinear Schlesinger
equations, i.e. the differential equations describing the
isomonodromy deformations of the Fuchsian systems.

In this paper we will study the Akhiezer polynomials $P_{n}$ using
the Riemann-Hilbert approach introduced in the theory of
orthogonal polynomials in \cite{FIK}. This will allow us
to exploit the well-developed Riemann-Hilbert and algebro-geometric
schemes of the Soliton theory \cite{nov}, \cite{fad}, \cite{bel} -
with certain important technical modifications
though, and not only re-derive the previous results of \cite{Chen}
and  \cite{mag} but also unite them in a single approach and produce
further facts concerning  the Akhiezer polynomials.
Specifically, in addition to the derivation of
Magnus's equations, we will solve them in terms of the multidimensional
$\Theta$-functions, and we will identify the corresponding Hankel determinant
with the relevant $\tau-$ function, i.e. with one of the
central objects associated with an integrable system, in
our case - with the Magnus-Schlesinger equation. It should also be
mentioned that part of our $\Theta$ - formulae, namely the ones describing
the recurrence coefficients
and the related Baker-Akhiezer function,
reproduce the known expressions obtained in the late 70s 
(the works of I. Krichever, D. Mumford, S. Novikov, and M. Salle) 
for the finite-gap discrete Schr\"odinger 
operators which were then intensively studied in connection with 
the periodic Toda lattice
(see the pioneering paper of H. Flaschka and D. McLaughlin \cite{flash} and
also  \cite{nov} and \cite{fad} for more 
on the history of the subject).
  
We would like to think of our paper as a tribute to
the pioneering works of N. I. Akhiezer which layed the foundation for
the construction, in the 1970's of the algebro-geometric method in the theory
of integrable systems, 
whose modern ``Riemann-Hilbert'' version we are using here.

\vskip .4cm
\setcounter{equation}{0}
\section{Riemann-Hilbert problem}

 According to the classical theory of orthogonal polynomials
the monic $P_n,$ (with $P_0=1$ and $P_{-1}=0$)  and
the polynomials of the second kind,
\bq
Q_n(z):=\int_{E}\frac{P_n(z)-P_n(t)}{z-t}w_+(t)dt,
\eq
of degree $n-1,$ are linearly independent solutions of the
second order difference equation,
\bq\label{rec}
zv_n(z)=v_{n+1}(z)+b_{n+1}v_n(z)+a_nv_{n-1}(z).
\eq
Following the general scheme of \cite{FIK} (see also \cite{BI}, \cite{DKMVZ}),
let us introduced the $2\times 2$ matrix
$Y_n(z)$ be defined for
$n=1,2,..$  and $z\in \Co$ as follows:
\bq \label{8}
\mbox{} Y_n(z)&=&\left(\matrix{P_n(z)&\int_{E}
\frac{P_n(t)\wp(t)}{z-t}dt\cr
\frac{P_{n-1}(z)}{h_{n-1}}&\frac{1}{h_{n-1}}\int_{E}
\frac{P_{n-1}(t)\wp(t)}{z-t}dt\cr}\right)\nonumber\\
\mbox{} &=&\left(\matrix{P_n(z)&\psi(z)P_n(z)-Q_n(z)\cr
\frac{P_{n-1}}{h_{n-1}}&\frac{\psi(z)P_{n-1}(z)-
                         Q_{n-1}(z)}{h_{n-1}}\cr}\right)
\eq
where
\bq
\mbox{} \psi(z)&:=&\int_{E}\frac{\wp(t)}{z-t}dt=
{\sqrt {\frac{\Pi_{i=1}^{g}(z-\al_i)}{\Pi_{j=0}^{g+1}(z-\bt_i)}}}\nonumber\\
\mbox{} &=&-i\pi w(z) .
\eq

{\bf Proposition 1.} The function  $Y_n(z)$ satisfies the
following conditions,
\bq
\mbox{} &{\bf {RH1.}}& Y_n(z)\;{\rm is\;analytic\;in\;}
\Co\setminus E\nonumber\\
\mbox{} &{\bf {RH2.}}& Y_{n,-}(z)=
Y_{n,+}(z)\left(\matrix{1&2\pi i\wp(z)\cr
                                        0&1\cr}\right),\;\;z\in E
\setminus \{\beta_{j}\}_{j=0}^{g+1}.
\nonumber\\
\mbox{} &{\bf {RH3.}}& Y_n(z)z^{-n\st}\to I,\;\;z\to\infty.\nonumber\\
\mbox{} &{\bf {RH4.}}&Y_n(z) = \hy^{(\bt_j)}_n(z)
\left(\matrix{\sqrt {z-\bt_j}&0\cr
              \frac{1}{\b_{j}}&\frac{1}{\sqrt {z-\bt_j}}\cr}\right),
\;\;z\;\in {\cal U}_{\bt_j},\;0\leq j\leq g+1,
\nonumber
\eq
where ${\cal U}_{z_{0}}$ denote a neighborhood of a point $z_{0}$.
The matrix valued function $\hy^{(\bt_j)}_n(z)$ is\\
\noindent
holomorphic in ${\sqrt {z-\bt_j}}$
and $\b_{j}$ is defined by the equation,
\be \label{bj}
w(z) = (z-\bt_j)^{-1/2}\b_{j} \frac{i}{\pi}(1 + O(z-\bt_{j})).
\ee
We shall also assume that the branch of ${\sqrt {z-\bt_j}}$
is defined by the condition,
$$
0 < \arg (z-\bt_{j}) < 2\pi, \quad \mbox{if}\mbox 
\quad j \leq g,\quad
\mbox{and}\quad
-\pi < \arg (z-\bt_{g+1}) < \pi.
$$
In addition, we assert that
\be \label{dethy}
\det \hy^{(\bt_j)}_n(\beta_{j}) = 1 \neq 0.
\ee
{\bf Proof.}

Using the basic properties of the Cauchy integrals and the Plemelj formulae
we directly verify that $Y_n(z)$ satisfies {\bf{RH1}} - {\bf{RH2}}. To check
property {\bf{RH3}} it is enough to note  that {\it because of the
orthogonality condition (\ref{ortog})}, we have (cf. \cite{FIK,Chen})
\bq
\int_{E}
\frac{P_n(t)\wp(t)}{z-t}dt&=& \sum_{k=0}^{\infty}\frac{1}{z^{k+1}}\int_{E}P_n(t)\wp(t)
t^{k}dt\nonumber\\
\mbox{} &=& \frac{h_n}{z^{n+1}} + {\rm O}\left(\frac{1}{z^{n+2}}\right),\quad z \to \infty.
\nonumber
\eq
To prove  {\bf{RH4}} we observe that the matrix product,
\bq
Y_{n}(z)\left(\matrix{\frac{1}{\sqrt {z-\bt_j}}&0\cr
             - \frac{1}{\b_{j}}&\sqrt {z-\bt_j}\cr}\right),\nonumber
\eq
is bounded near $\bt_j$ (the singular terms in the first column
cancel out), and hence the function $\hy^{(\bt_j)}_n(z)$ defined
by equation {\bf{RH4}} is indeed holomorphic in ${\sqrt {z-\bt_j}}$.
To complete the prove of the proposition we only need to establish
equation (\ref{dethy}). To this end, we notice that we have already established
{\bf{RH1}} - {\bf{RH4}} but short of equation (\ref{dethy}). One can see, however,
that {\bf{RH1}} - {\bf{RH4}} already yield even stronger statement. Namely, we
claim that
\be \label{dety}
\det Y_{n}(z) \equiv 1.
\ee
Indeed, the (scalar) function $\det Y_{n}(z)$ is holomorphic in
$\cpo\setminus E$,
has no jumps across $E$ and has removable singularities at the end
points of $E$; moreover, it
approaches $1$ as $z\longrightarrow\infty$. By the Liouville theorem,
equation (\ref{dety}) follows. Equation (\ref{dethy}) is a direct
consequence of equation  (\ref{dety}). The proposition is proven.

{\bf Remark 2.1} Equation (\ref{dety}) can be also derived by using
the first line of (\ref{8}) and the Christoeffel-Darbooux formula,
\bq
\mbox{} \det Y_n(z)&=&\frac{1}{h_{n-1}}
\int_{E}\frac{P_n(z)P_{n-1}(t)-P_{n-1}(z)P_n(t)}{z-t}\wp(t)dt
\nonumber\\
\mbox{} &=&\int_{E}\sum_{k=0}^{n-1}\frac{1}{h_k}P_k(z)P_k(t)
\wp(t)dt\nonumber\\
\mbox{} &=&\int_{E}K_{n}(z,t)\wp(t)dt=P_0(z)h_0=1,
\eq
or from the recurrence relations,
\be
\det Y_n(z)=\frac{1}{h_{n-1}}
(Q_n(z)P_{n-1}(z)-P_n(z)Q_{n-1}(z))=1.
\nonumber
\ee

{\bf Remark 2.2} $Y_n(z)$ also depend on
$\{\de_j:1\leq j\leq 2g+2\}.$

\vskip .2in

{\bf Proposition 2.} Conditions {\bf{RH1}} - {\bf{RH4}} defines the function
$Y_{n}(z)$ uniquely.

{\bf Proof.} If $\ty_n(z)$ is another function
that satisfies {\bf{RH1}} - {\bf{RH4}} then
$X_n(z):={\tilde Y}_n(z)Y_n^{-1}(z)$ is
holomorphic for $z\in\cpo\setminus\{\bt_j:0\leq j\leq g+1\}.$
Furthermore, for $z\;\in {\cal U}_{\bt_j},$
\be
Y_n^{-1}(z) = \left(\matrix{\frac{1}{\sqrt {z-\bt_j}}&0\cr
              -\frac{1}{\b_{j}}&\sqrt {z-\bt_j}\cr}\right)
\hy^{(\bt_j)-1}_n(z),
\ee
where $\hy^{(\bt_j)-1}_n(z)$ is holomorphic (see equation (\ref{dethy}) !)
in ${\sqrt {z-\bt_j}}.$ This implies,
\be
X_n(z)= O(1), \quad z \sim \bt_{j},
\ee
which in turn implies $X_n(z)$ is holomorphic for
$ z\in\cpo,$ and $X_n(z)=I,$ for all $z\in\cpo.$

\vskip .2in
The conditions {\bf{RH1}} - {\bf{RH4}} constitute the Riemann-Hilbert
problem whose unique solution is given by equation (2.3), due to Proposition 1.
\vskip .2in

The Riemann-Hilbert problem {\bf{RH1}} - {\bf{RH4}} together with the equation
\be\label{P_Y}
P_{n}(z) = (Y_{n}(z))_{11}
\ee
will be used as an {\it alternative definition} of the Akhiezer polynomials.
Notice also that the asymptotic  condition {\bf{RH3}} can be extended to the full
Laurent series,
\bq\label{Yinf}
Y_n(z)=\left(I+\sum_{k=1}^{\infty}\frac{m_k(n)}{z^k}\right)z^{n\st},
\quad |z| > 1
\eq
and from (2.3) we have,
\bq
m_1(n)=\left(\matrix{p_1(n)&h_n\cr
                     1/h_{n-1}&-p_1(n)\cr}\right)
\eq
where $p_1(n)$ is the coefficient of $z^{n-1}$ of $P_n(z).$
Taking into account the recurrence relations (\ref{rec}), we have,
$$
a_{n} = \frac{h_{n}}{h_{n-1}}
$$
and
$$
b_{n+1} = p_{1}(n) - p_{1}(n+1),
$$
and the following relations supplementing (\ref{P_Y})
\bq\label{hab_Y}
\mbox{} h_n&=&\left(m_{1}(n)\right)_{12}\label{hab_Y1}\\
\mbox{} a_n&=&\left(m_{1}(n)\right)_{12}\left(m_{1}(n)\right)_{21}\label{hab_Y2}\\
\mbox{} b_{n+1}&=&\left(m_{1}(n)\right)_{11}
-\left(m_{1}(n+1)\right)_{11}.\label{hab_Y3}
\eq
Therefore, all the basic ingredients of the theory of polynomials $P_{n}(z)$
(including the polynomials themselves)
can be obtained directly from the solution $Y_{n}(z)$ of the
Riemann-Hilbert problem.

{\bf Remark 2.3} In the {\it a prior} setting of the Riemann-Hilbert
problem {\bf{RH1}} - {\bf{RH4}}, the condition {\bf{RH3}} can be
replaced by the following weaker one
\bq
\mbox{} &{\bf {RH4.}}&Y_n(z)
\left(\matrix{\frac{1}{\sqrt {z-\bt_j}}&0\cr
              -\frac{1}{\b_{j}}&\sqrt {z-\bt_j}\cr}\right)
 = O(1), \; z \sim \bt_{j},\;0\leq j\leq g+1.
\eq

\setcounter{equation}{0}
\section{Differential Equations}

Having obtained equations (\ref{P_Y}) - (\ref{hab_Y3}) which represent
orthogonal polynomials $P_{n}(z)$ and the corresponding norm and
recurrence coefficients in terms of the solution $Y_{n}(z)$
of the Riemann-Hilbert problem {\bf{RH1}} - {\bf{RH4}}, we
can now use the powerfull techniques of the Soliton theory.
Specifically, in this and the two following sections we will apply
a certain modification of the standard Zakharov-Shabat dressing
method (see e.g. \cite{nov}) to obtain the relevant differential
and difference equations for the Akhiezer polynomials. The modification
needed is caused by the presence of the condition {\bf{RH4}}.
This condition indicates the relation
of the problem under consideration to the theory of Fuchsian
systems. Indeed, our derivations will be close to
the Zakharov - Shabat scheme and to the constructions of the Jimbo-Miwa-Ueno
monodromy theory \cite{jmu} (see also \cite{its1} were both
methods are unified in a single general Riemann-Hilbert formalism).

 To describe the change of $Y_n(z)$ with respect to $z$ for a
fixed $n,$ it is advantageous to transform the Riemann-Hilbert
problem satisfied by $Y_n(z)$ in to a form where jump matrix
has constant entries. To this end, put
\be \label{Phidef}
\Phi_n(z)=Y_n(z)\left(\matrix{1&0\cr
                          0&w^{-1}(z)\cr}\right)
            \left(\matrix{\stp&0\cr
                          0&1/\stp\cr}\right).
\ee
A direct computation shows that
\bq\label{Phi_jump}
\mbox{} \Phi_{n,-}(z)=
\Phi_{n,+}(z)
\left(\matrix{1&-1\cr
              0&-1\cr}\right).
\eq
To specify the behavior of the new function near the end points
of the set $E$ let us observe that the new (constant !) jump matrix
admits the following spectral representation,
$$
\left(\matrix{1&-1\cr
              0&-1\cr}\right) = {\bf P}^{-1}
\left(\matrix{-1&0\cr
              0&1\cr}\right){\bf P},
$$
where
$$
{\bf P} = \left(\matrix{0&1\cr
              2&-1\cr}\right).
$$
This implies that the function
$$
\Phi^{(\bt_{j})}(z) := \left(\matrix{\sqrt{z -\bt_{j}}&0\cr
              0&1\cr}\right){\bf P}
$$
satisfies the jump condition (\ref{Phi_jump}) in the neighborhood of $\bt_{j}$. Indeed,
assuming $z \in {\cal U}_{\bt_{j}}\cap E,$ we find,
$$
[\Phi_{+}^{(\bt_{j})}(z)]^{-1}\Phi_{-}^{(\bt_{j})}(z)
={\bf P}^{-1}\left(\matrix{\frac{1}{\left(\sqrt{z -\bt_{j}}\right)_{+}}&0\cr
              0&1\cr}\right)
\left(\matrix{\left(\sqrt{z -\bt_{j}}\right)_{-}&0\cr
              0&1\cr}\right){\bf P}
$$
$$
={\bf P}^{-1}
\left(\matrix{-1&0\cr
              0&1\cr}\right){\bf P}
=\left(\matrix{1&-1\cr
              0&-1\cr}\right).
$$
Hence the matrix valued fucntion
$$
\Phi_{n}(z)[\Phi^{(\bt_{j})}(z)]^{-1}
$$
has no jump accross $E$ and therefore is holomorphic in the punctured neighborhood
${\cal U}_{\bt_{j}}\setminus \{\bt_{j}\}$. Observe in addition that in the product,
\be
\left(\matrix{\sqrt{z-\bt_{j}}&0\cr
            \frac{2\pi i}{\b_{j}}&\frac{1}{\sqrt{z-\bt_{j}}}\cr}\right)
\left(\matrix{1&0\cr
            0&w^{-1}(z)\cr}\right)
\left(\matrix{\frac{1}{2}&\frac{1}{2}\cr
            1&0\cr}\right)
\left(\matrix{\frac{1}{\sqrt{z-\bt_{j}}}&0\cr
            0&1\cr}\right),
\ee
the negative powers of  $\sqrt{z-\bt_{j}}$ cancel out.

 Therefore we conclude that the product
$\Phi_{n}(z)[\Phi^{(\bt_{j})}(z)]^{-1}$ is in fact holomorphic in the
whole neighborhood  $ {\cal U}_{\bt_{j}}$.
Similar is also true for the matrix product
$$
\Phi_{n}(z)[\Phi^{(\al_{j})}(z)]^{-1} \equiv
\Phi_{n}(z)\left[
\left(\matrix{\frac{1}{\sqrt{z -\al_{j}}}&0\cr
              0&1\cr}\right){\bf P}
\right]^{-1}
$$
in the neighborhood ${\cal U}_{\al_{j}}$ of the endpoint $\al_{j}$.
Here we shall assume that the branch of $\sqrt{z -\al_{j}}$ is defined
by the condition,
$$
-\pi < \arg(z - \al_{j}) < \pi.
$$
In summary,  $\Phi_n(z)$ solves the following
Riemann-Hilbert problem:
\bq
\mbox{} &{\bf \Phi 1}.&\; \Phi_n(z)\;{\rm is\;holomorphic\;for\;}
z\in\Co \setminus E\nonumber\\
\mbox{} &{\bf \Phi 2}.&\; \Phi_{n,-}(z)=\Phi_{n,+}(z)
\left(\matrix{1&-1\cr
              0&-1\cr}\right),\;z\in E
\nonumber\\
\mbox{} &{\bf \Phi 3}.&\; \Phi_n(z)=
\left(I+{\rm O}\left(\frac{1}{z}\right)\right)
                 z^{\left(\matrix{n&0\cr
                                  0&-n+1\cr}\right)}
\left(\matrix{\sqrt{2\pi i}&0\cr
              0&-\sqrt{\frac{\pi i}{2}}\cr}\right),\;\;
z\longrightarrow \infty, \nonumber\\
\mbox{} &{\bf \Phi 4}.&\;
\Phi_n(z) = \hat{\Phi}^{(\bt_{j})}_{n}(z)
\left(\matrix{\sqrt{z-\bt_{j}}&0\cr
            0&1\cr}\right)
\left(\matrix{0&1\cr
            2&-1\cr}\right),\;\;\quad z\in {\cal U}_{\bt_{j}}\nonumber\\
\mbox{} &{\bf \Phi 5}.&\;
\Phi_n(z) = \hat{\Phi}^{(\al_{j})}_{n}(z)
\left(\matrix{\frac{1}{\sqrt{z-\al_{j}}}&0\cr
            0&1\cr}\right)
\left(\matrix{0&1\cr
            2&-1\cr}\right),\;\;\;
\quad z\in {\cal U}_{\al_{j}},\nonumber
\eq
where $\hat{\Phi}^{(\bt_{j})}_{n}(z)$ and $\hat{\Phi}^{(\al_{j})}_{n}(z)$
are holomorphic in the neighborhoods of the points
$\bt_{j}$ and $\al_{j}$, respectively.
Moreover, the matrices $\hat{\Phi}^{(\bt_{j})}_{n}(\bt_{j})$ and
$\hat{\Phi}^{(\al_{j})}_{n}(\al_{j})$ are invertible. In fact,
$$
\hat{\Phi}^{(\bt_{j})}_{n}(\bt_{j}) = \hat{Y}^{(\bt_j)}_{n}(\beta_{j})
\left(\matrix{\sqrt{\frac{\pi i}{2}}&0\cr
            0&\frac{1}{\b_{j}}\sqrt{\frac{\pi i}{2}}\cr}\right)
$$
and
$$
\hat{\Phi}^{(\al_{j})}_{n}(\al_{j}) = Y_{n}(\al_{j})
\left(\matrix{0&\sqrt{\frac{\pi i}{2}}\cr
            -\frac{1}{\a_{j}}\sqrt{\frac{\pi i}{2}}&0\cr}\right)
$$
where $\a_{j}$ is defined by the equation (cf. \ref{bj})
\be \label{aj}
w(z) = (z-\al_j)^{1/2}\a_{j} \frac{i}{\pi}(1 + O(z-\al_{j})).
\ee
We want to emphasize, that unlike the case of the $Y$ - Riemann-Hilbert problem,
in the case of the $\Phi$ - Riemann-Hilbert problem the left multipliers
$\hat{\Phi}^{(\bt_{j})}_{n}(z)$
and $\hat{\Phi}^{(\al_{j})}_{n}(z)$ are holomorphic {\it with respect to $z$}.

{\bf Remark 3.1}
 From ${\bf \Phi 1 - \Phi 5}$ it follows
(independent of (\ref{Phidef})) that
\bq
\det\Phi_n(z)=\frac{1}{w(z)}.
\eq

\vskip .2in

Consider now, the logarithmic derivative of $\Phi_n(z),$
\bq
A(z,n)&:=&\frac{d\Phi_n(z)}{dz}\Phi_n^{-1}(z).
\eq
Since all the right matrix multipliers in the r.h.s. of
${\bf \Phi 2 - \Phi 5}$ are constant matrices, $A(z,n)$ enjoys the following
properties:
\bq
\mbox{} &A1.&\;A(z,n){\;\rm is\;}{\rm holomorphic\;for\;}z
\in\cpo\setminus\{\al_j,\bt_j\},\nonumber\\
\mbox{} &A2.&\;A(z,n)=
\frac{\left(\matrix{n&0\cr
                    0&-n+1\cr}\right)}{z}+
{\rm O}\left(\frac{I}{z^2}\right),\;\;
z\longrightarrow\infty,\nonumber\\
\mbox{} &A3.&\;A(z,n)= \frac{1}{2}\hat{\Phi}^{(\bt_{j})}_{n}(\bt_{j})
\frac{\left(\matrix{1&0\cr
            0&0\cr}\right)}{z-\bt_{j}}\hat{\Phi}^{(\bt_{j})-1}_{n}(\bt_{j})
+ {\rm O}(1), \quad z\sim \bt_{j}, \nonumber\\
\mbox{} &A4.&\;A(z,n)= - \frac{1}{2}\hat{\Phi}^{(\al_{j})}_{n}(\al_{j})
\frac{\left(\matrix{1&0\cr
            0&0\cr}\right)}{z-\al_{j}}\hat{\Phi}^{(\al_{j})-1}_{n}(\al_{j})
+ {\rm O}(1), \quad z\sim \al_{j}.\nonumber
\eq
By virtue of the Liouville theorem, it follows that,
\be
A(z,n)=\sum_{j=0}^{g+1}\frac{B_j(n)}{z-\bt_j}+
\sum_{j=1}^{g}\frac{A_j(n)}{z-\al_j}
\ee
where
\be \label{Bj}
B_{j}(n) := \frac{1}{2}\hat{\Phi}^{(\bt_{j})}_{n}(\bt_{j})
\left(\matrix{1&0\cr
            0&0\cr}\right)\hat{\Phi}^{(\bt_{j})-1}_{n}(\bt_{j})
=\frac{1}{2}\hat{Y}^{(\bt_j)}_{n}(\bt_{j})
\left(\matrix{1&0\cr
            0&0\cr}\right)\hat{Y}^{(\bt_j)-1}_{n}(\bt_{j})
\ee

\be \label{Aj}
A_{j}(n) := -\frac{1}{2}\hat{\Phi}^{(\al_{j})}_{n}(\al_{j})
\left(\matrix{1&0\cr
            0&0\cr}\right)\hat{\Phi}^{(\al_{j})-1}_{n}(\al_{j})
=-\frac{1}{2}Y_{n}(\al_{j})
\left(\matrix{1&0\cr
            0&0\cr}\right)Y^{-1}_{n}(\al_{j}).
\ee
Note also,
\bq
\mbox{} \sum_{j=0}^{g+1}B_j(n)+\sum_{j=1}^{g}A_j(n)=
\left(\matrix{n&0\cr
              0&-n+1\cr}\right).\nonumber
\eq
Using (2.3) and {\bf {RH4}} give
\bq
\hy^{(\bt_j)}_n(\bt_j)=\left(\matrix{Q_n(\bt_j)/\b_j&\b_jP_n(\bt_j)\cr
            \frac{Q_{n-1}(\bt_j)}{\b_jh_{n-1}}&\b_jP_{n-1}(\bt_j)/h_{n-1}}
\right).
\eq
 We conclude this section by recording the linear matrix
differential equation with Fuchsian singularities at
$\{\al_j,\bt_j\},$ mentioned in the Abstract,
\bq\label{fuchs}
\mbox{} \frac{d\Phi_n(z)}{dz}&=&A(z,n)\Phi_n(z),
\eq
with $A(z,n)$ defined by (3.7), (3.8) and (3.9).
  Furthermore, using the second line of (2.3), the matrix valued residues are
expressed in terms of the evalutions of the polynomials at the branch
points:
\bq
\mbox{} B_j(n)&=&\frac{1}{2}
\left(\matrix{Q_n(\bt_j)P_{n-1}(\bt_j)/h_{n-1}&-Q_n(\bt_j)P_n(\bt_j)\cr
 Q_{n-1}(\bt_j)P_{n-1}(\bt_j)/h_{n-1}^2&
   -Q_{n-1}(\bt_j)P_n(\bt_j)/h_{n-1}\cr}
\right)\label{BjPQ}\\\nonumber\\
\mbox{} A_j(n)&=&\frac{1}{2}
\left(\matrix{P_n(\al_j)Q_{n-1}(\al_j)/h_{n-1}&-Q_n(\al_j)P_n(\al_j)\cr
Q_{n-1}(\al_j)P_{n-1}(\al_j)/h_{n-1}^2&
-P_{n-1}(\al_j)Q_{n}(\al_j)/h_{n-1}\cr}\right).\label{AjPQ}
\eq
Note that from (\ref{Bj}) and (\ref{Aj}) it follows that
\bq
\mbox{} {\rm tr}B_j(n)&\equiv&\frac{1}{2 h_{n-1}}
(Q_{n-1}(\bt_j)P_n(\bt_j)-Q_n(\bt_j)P_{n-1}(\bt_j))=1/2\nonumber\\
\mbox{} \det B_j(n)&=&0,\nonumber
\eq
and
\bq
\mbox{} {\rm tr}A_j(n)&\equiv&-\frac{1}{2h_{n-1}}
(Q_n(\al_j)P_{n-1}(\al_j)-P_n(\al_j)Q_{n-1}(\al_j))=-1/2.\nonumber\\
\mbox{} \det A_j(n)&=&0.\nonumber
\eq
We note that this leads to a discrete analogue of the ``Wronskian'' relation,
\bq
P_{n-1}(z)Q_n(z)-P_{n}(z)Q_{n-1}(z)=h_{n-1},\nonumber
\eq
which, of course, can be independently derived from the recurrence relations.

 As it has already been mentioned in Introduction, equation
(\ref{fuchs}), even for more general weights of the type $\prod_j(t-\de_j)^{\kappa_j}$,
was first obtained in  \cite{mag}. In \cite{mag} the Riemann-Hilbert
problem is not used explicitely; rather, the author analyses directly
the monodromy properties of the function $Y_{n}(z)$, i.e. the approach
of \cite{mag} is based more on the ideas of \cite{jmu} than of \cite{nov}.
It is also worth mentioning that our approach can be extended to the
general semi-classical weights without any serious modifications.

\setcounter{equation}{0}
\section{{\bf Derivatives with respect to the branch points.}}

 In this section we determine differentiation formulas for
$\Phi_n(z)$ with respect to $\{\al_j,\bt_j\}.$ First let us
consider the logarithmic derivative of $\Phi_n(z)$ with respect to a
particular $\bt_j;$
\bq
\mbox{} V_j(z):=\frac{\pl\Phi_n(z)}{\pl\bt_j}\Phi^{-1}_n(z),
\eq
and note that $V_j(z)$ has the following properties
\bq
\mbox{} &V1.&\;V_j(z)\;
{\rm is\;holomorphic\;for\;}z\in {\bf C}\setminus
\{\bt_j\}.
\nonumber\\
\mbox{} &V2.&\;V_j(z)={\rm O}(I/z),\quad z\to\infty.\nonumber\\
\mbox{} &V3.&\;V_j(z)\sim
-\frac{1}{2}\hat{\Phi}^{(\bt_{j})}_{n}(\bt_{j})
\frac{\left(\matrix{0&0\cr
                    0&1\cr}\right)}{z-\bt_j}\hat{\Phi}^{(\bt_{j})-1}_{n}
(\bt_{j})+ {\rm O}(1).\nonumber\\
\eq
By comparing with (\ref{Bj}) and again invoking
the Liouville theorem, we conclude that
\bq
V_j(z)=-\frac{B_j(n)}{z-\bt_j},
\eq
which implies
\bq
\pl_{\bt_j}\Phi_n(z)&=&-\frac{B_j(n)}{z-\bt_j}\Phi_n(z).
\eq
A similar analysis gives
\bq
\pl_{\al_j}\Phi_n(z)&=&-\frac{A_j(n)}{z-\al_j}\Phi_n(z)
\eq
\setcounter{equation}{0}
\section{{\bf Difference Equation.}}
Consider the ``difference logarithmic derivative''
$$
U_n(z):= \Phi_{n+1}(z)\Phi^{-1}_{n}(z) \equiv Y_{n+1}(z)Y^{-1}_{n}(z).
$$
Taking into account that all the right matrix multipliers
in the r.h.s of {\bf{RH1}} - {\bf{RH4}} are constant with respect
to $n$ we conclude that $U_n(z)$ is an entire function.
Moreover, from the asymptotics (\ref{Yinf}) we have that
\bq
\mbox{} U_n(z)&=&\left( I + \frac{m_{1}(n+1)}{z}\right)z^{\sigma_{3}}
\left( I - \frac{m_{1}(n)}{z}\right) + O\left(\frac{1}{z}\right)\nonumber\\
\mbox{} &=&\left( I + \frac{m_{1}(n+1)}{z}\right)
\left(\matrix{z&0\cr
              0&0\cr}\right)
\left( I - \frac{m_{1}(n)}{z}\right) + O\left(\frac{1}{z}\right)
\nonumber\\
\mbox{}&=&z\left(\matrix{1&0\cr
                       0&0\cr}\right)
+ m_{1}(n+1)\left(\matrix{1&0\cr
                       0&0\cr}\right)
- \left(\matrix{1&0\cr
                       0&0\cr}\right)m_{1}(n)
+ O\left(\frac{1}{z}\right),\quad z \to \infty.\nonumber
\eq
Appealing once again to the Liouville theorem, we conclude that
$U_n(z)$ is linear function in $z$ defined by the equations
\bq
\mbox{} U_n(z)&=&\left(\matrix{z+(m_{1}(n+1))_{11} - (m_{1}(n))_{11}&- (m_{1}(n))_{12}\cr
                (m_{1}(n+1))_{21}&0\cr}\right)\nonumber\\
\mbox{} &=&\left(\matrix{z-b_{n+1}&-h_n\cr
                       1/h_n&0\cr}\right),\nonumber
\eq
where in the last equation we have taken into account (\ref{hab_Y1})- (\ref{hab_Y3}).
To summarize, the difference equation for the function $\Phi_n(z)$ reads
\bq\label{difference}
\Phi_{n+1}(z)=\left(\matrix{z-b_{n+1}&-h_n\cr
                       1/h_n&0\cr}\right)\Phi_n(z).
\eq

Of course, equation (\ref{difference}) is just the matrix form of the basic
recurrence equation (\ref{rec}). Nevertheless, we gave its ``Riemann-Hilbert''
derivation to emphasize the ``master'' role of the Riemann-Hilbert problem
{\bf{RH1}} - {\bf{RH4}} in our analysis.

\setcounter{equation}{0}
{\section {Schlesinger Equations and the Hankel Determinant.}}

 With the unified notation mentioned in the Introduction,
we write,
\be
A(z,n)=\sum_{j=1}^{2g+2}\frac{C_j(n)}{z-\de_j},
\ee
and the correspondence,
\be
(A_1(n),...,A_g(n),B_{0}(n),...,B_{g+1}(n))
\longrightarrow(C_1(n),...,C_g(n),
C_{g+1}(n),...,C_{2g+2}(n)).
\ee
Note that, $C_j(n),$ depend on $\de_j'$s.
We of course have,
\bq
\mbox{}\pz\Phi_n(z)&=&
\sum_{j=1}^{2g+2}\frac{C_j(n)}{z-\de_j}\Phi_n(z),\\
\mbox{}\pdk\Phi_n(z)&=&-\frac{C_k(n)}{z-\de_k}\Phi_n(z).
\eq
Applying $\pz$ on (6.4) gives
\bq
\pz\pdk\Phi_n(z)=\frac{C_k(n)}{(z-\de_k)^2}\Phi_n-
\frac{C_k(n)}{z-\de_k}\sum_{j=1}^{2g+2}\frac{C_j(n)}{z-\de_j}
\Phi_n,
\eq
and $\pdk$ on (6.3) gives,
\bq
\pdk\pz\Phi_n(z)=\frac{C_k(n)}{(z-\de_k)^2}\Phi_n+
\sum_{j=1}^{2g+2}\frac{\pdk C_j(n)}{z-\de_j}\Phi_n
-\left(\sum_{j=1}^{2g+2}\frac{C_j(n)}{z-\de_j}\right)
\frac{C_k(n)}{z-\de_k}\Phi_n.
\eq
Since $\pz\pdk\Phi_n=\pdk\pz\Phi_n$ and $\det\Phi_n\neq 0,$
we get,
\bq
\mbox{}\sum_{j=1}^{2g+2}\frac{\pdk C_j(n)}{z-\de_k}
&=&\sum_{j=1}^{2g+2}
\frac{\left[C_j(n),C_k(n)\right]}{(z-\de_j)(z-\de_k)}
\nonumber\\
\mbox{}&=&\sum_{j=1}^{2g+2}
\frac{\left[C_j(n),C_k(n)\right]}{\de_j-\de_k}
\left(\frac{1}{z-\de_j}-\frac{1}{z-\de_k}\right).
\eq
We now send $z$ to a particular $\de_j$ in (6.7), with $j\neq k,$
and find by equating residues,
\bq
\pdk C_j(n)=\frac{\left[C_j(n),C_k(n)\right]}{\de_j-\de_k},
\;\;\;j\neq k.
\eq
If $j=k$, then a similar calculation gives,
\bq
\pdk C_k(n)=-\sum_{l(\neq k)}
\frac{\left[C_l(n),C_k(n)\right]}{\de_l-\de_k}.
\eq
The equations (6.8) and (6.9) are the
Schlesinger Equations satisfied by $C_j(n)$.
This is the equation first derived for the general semi-classical
orthogonal polynomials in \cite{mag}. We are now going to move
beyond the results of \cite{mag} and show that the corresponding
$\tau$ - function can be identified with the Hankel determinant
associated with the weight $w_{+}(t)$. To this end we first
recall Jimbo-Miwa-Ueno definition of the $\tau$ -function.

Let $\Omega^{(1)}$ be the one-form,
\bq
\mbox{} \Omega^{(1)}(\de_1,...,\de_{2g+2})
&:=&\sum_{1\leq j<k\leq 2g+2}{\rm tr}\left(C_j(n)C_k(n)\right)
\frac{d\de_j-d\de_k}{\de_j-\de_k}\nonumber\\
\mbox{} &=&\sum_{1\leq j<k\leq 2g+2}{\rm tr}\left(C_j(n)C_k(n)\right)
d\ln|\de_j-\de_k|,
\eq
then it can be verified \cite{jmu} using the Schlesinger Equations
that,
\bq
d\Omega^{(1)}=0,
\eq
which implies that, localy, $\Omega^{(1)}$ is an exact form.
The  $\tau-$ function of the completely
integrable system of partial differential equations (6.8) and (6.9)
is then defined by the relation,
\bq
\Omega^{(1)}=d\ln\tau_n(\de_1,...,\de_{2g+2}).
\eq

In the theory orthogonal
polynomials, the Hankel determinant,
\bq
D_n[w_+]:=\det\left(\int_{E}t^{j+k}\wp(t)dt\right)_{j,k=0}^{n-1},
\eq
has two other equivalent expressions,
\bq
\mbox{} D_n[w_+]&=&\frac{1}{n!}\int_{E}...\int_E
\prod_{1\leq j<k\leq n}(x_j-x_k)^2\prod_{l=1}^{n}\wp(x_l)dx_l,
\nonumber\\
\mbox{} &=&\prod_{j=0}^{n-1}h_j.
\eq
 It is to be expected from the structure of
the Riemann-Hilbert formulation that, $D_n,$ considered as
a function of $\{\de_j\}_{j=1}^{2g+2},$ is the $\tau-$
function for this problem. To understand this,
we require the derivatives of $h_n$ w.r.t. to $\de_k.$
To begin with, we use that,
\be
\pdk\Phi_n(z)=-\frac{C_k(n)}{z-\de_k}\Phi_n(z),
\ee
must be compatible with
\be
\Phi_{n+1}(z)=\left(\matrix{z-b_{n+1}&-h_n\cr
                          1/h_n&0\cr}\right)\Phi_n(z).
\ee
This results is
\bq
\mbox{}\left(\matrix{z-b_{n+1}&-h_n\cr
      h_n&0\cr}\right)\frac{C_k(n)}{z-\de_k}
&-&\frac{C_k(n+1)}{z-\de_k}\left(\matrix{z-b_{n+1}&-h_n\cr
                                 h_n&0\cr}\right)\nonumber\\
&=&\left(\matrix{-\pdk b_{n+1}&-\pdk h_n\cr
               -(1/h_n)\pdk\ln h_n&0\cr}\right),
\eq
which holds for all $z \in \cpo \setminus \{\de_{1}, ...
\de_{2g+2}\}$. Putting $z=\infty$ in (6.17), gives,
\bq
\left(\matrix{C_k^{11}(n)-C_k^{11}(n+1)&C_k^{12}(n)\cr
              -C_k^{21}(n+1)&0\cr}\right)
=\left(\matrix{-\pdk b_{n+1}&-\pdk h_n\cr
               -(1/h_n)\pdk\ln h_n&0\cr}\right),\nonumber
\eq
which implies, amongst others,
\bq
\pdk h_n=-C_k^{12}(n).
\eq
\vskip .4cm
\noindent
{\bf Lemma 1.} Let the asymptotic expansion of
$A(z,n)$ about $z=\infty$ be
\bq
A\left(z,n\right)=\sum_{k=0}^{\infty}\ca_k(n)z^{-k-1}
\eq
where
\bq
\ca_k(n):=\sum_{j=1}^{2g+2}C_j(n)\de^k_j(n).
\eq
Then the first two $\ca_k(n)$ are
\bq
\mbox{}\ca_0(n)&=&\left(\matrix{n&0\cr
                                0&1-n\cr}\right),\\
\mbox{}\ca_1(n)&=&\left(\matrix{0&0\cr
                                0&c_1\cr}\right)
+m_1(n)\left(\matrix{n-1&0\cr
                     0&-n\cr}\right)
-\left(\matrix{n&0\cr
               0&1-n\cr}\right)m_1(n),
\eq
where $c_{1} = \sum_{j=1}^{g}\left(\beta_{j} - \alpha_{j}\right)$.

{\bf Proof.} Putting (3.1) into (3.11) we find
\bq
\frac{d}{dz}Y_n(z)+Y_n(z)
\left(\matrix{0&0\cr
              0&-\frac{d}{dz}\ln w(z)\cr}\right)
=A(z,n)Y_n(z).
\eq
Expansion of (6.23) in $z^{-1}$ gives the desired results.
\vskip .4cm

\noindent
{\bf Theorem 2.} The Hankel determinant is the $\tau$ function
of the Magnus - Schlesinger Equations.

Proof: We start by equating the residues of (6.17)
at $z=\de_j$. This gives,
\bq
U_n(\de_j)C_j(n)=C_j(n+1)U_n(\de_j),\nonumber
\eq
or
\bq\label{freud1}
C_j(n+1)=U_n(\de_j)C_j(n)U_n^{-1}(\de_j),
\eq
where
\bq
\mbox{} U_n(z)&:=&\left(\matrix{z-b_{n+1}&-h_n\cr
                              1/h_n&0\cr}\right)\nonumber\\
\mbox{} U_n^{-1}(z)&=&\left(\matrix{0&h_n\cr
                                    -1/h_n&z-b_{n+1}\cr}\right).
\eq
A simple calculation shows that
\bq
\mbox{}U_n^{-1}(z)U_n(z^{\prime})=
\left(\matrix{1&0\cr
  \frac{z-z^{\prime}}{h_n}&1\cr}\right)
=I+\frac{z-z^{\prime}}{h_n}\s_{-},
\eq
where $\s_{-}:=\left(\matrix{0&0\cr
                            1&0\cr}\right).$
Now,
\bq
d\ln\tau_n=\sum_{k}\pdk\ln\tau_n\;\:d\de_k,
\eq
where (from (6.10)),
\bq
\pdj\ln\tau_n=\sum_{k(\neq j)}
\frac{\tr C_j(n)C_k(n)}{\de_j-\de_k},
\eq
which leads to
\bq
\mbox{}\pdj\ln\frac{\tau_{n+1}}{\tau_n}
&=&\sum_{k(\neq j)}\frac{\tr(C_j(n+1)C_k(n+1)-C_j(n)C_k(n))}
{\de_j-\de_k}\nonumber\\
&=&\sum_{k(\neq j)}
\frac{\tr(U_n(\de_j)C_j(n)U_n^{-1}(\de_j)U_n(\de_k)C_k(n)U_n^{-1}
(\de_k)-C_j(n)C_k(n))}{\de_j-\de_k}\nonumber\\
&=&\sum_{k(\neq j)}\frac{\tr(U_n^{-1}(\de_k)U_n(\de_j)C_j(n)
U_n^{-1}(\de_j)U_n(\de_k)C_k(n)-C_j(n)C_k(n))}{\de_j-\de_k}
\nonumber\\
&=&\sum_{k(\neq j)}\frac{\tr[(I-\frac{\de_j-\de_k}{h_n}
\s_{-})C_j(n)(I+\frac{\de_j-\de_k}{h_n}\s_{-})C_k(n)-
C_j(n)C_k(n)]}{\de_j-\de_k}.
\eq
A calculation shows that the term $[...]$ in  (6.29) is
\bq
\frac{\de_j-\de_k}{h_n}\left(C_j(n)\sm C_k(n)
-\sm C_j(n)C_k(n)\right)-\left(\frac{\de_j-\de_k}{h_n}\right)^2
\sm C_j(n)\sm C_k(n).\nonumber
\eq
We also note here some useful identities;
\bq
\mbox{}\tr(C_j(n)\sm C_k(n)&-&\sm C_j(n)C_k(n))\nonumber\\
\mbox{} &=&C_j^{12}(n)(C_k^{11}(n)-C_k^{22}(n))-C_k^{12}(n)
(C_j^{11}(n)-C_j^{22}(n)),\nonumber
\eq
and
\bq
\tr(\sm C_j(n)\sm C_k(n))=C_j^{12}(n)C_k^{12}(n).\nonumber
\eq
Therefore
\bq
\mbox{}\pdj\ln\frac{\tau_{n+1}}{\tau_n}
&=&\frac{1}{h_n}\sum_{k(\neq j)}\left(C_j^{12}(n)
(C_k^{11}(n)-C_k^{22}(n))-C_k^{12}(n)(C_j^{11}(n)-C_j^{22}(n))
\right)\nonumber\\
&-&\frac{1}{h_n^2}\sum_{k(\neq j)}(\de_j-\de_k)C_j^{12}(n)
C_k^{12}(n).
\eq
To simplify the r.h.s. of (6.30) we note, from (6.20), (6.21) and (6.22)
\bq
\mbox{} \sum_{j}C^{12}_j(n)&=&0,\nonumber\\
\sum_{j}(C_j^{11}(n)-C_j^{22}(n))&=&2n-1,\nonumber\\
\sum_{j}\de_jC_j^{12}(n)&=&-2nh_n.\nonumber
\eq
Using these, and $\sum_{k(\neq j)}f_k=-f_j+\sum_{k}f_k,$
the r.h.s. of (6.30), becomes,
\bq
\frac{C_j^{12}(n)}{h_n}\sum_k(C_k^{11}(n)-C_k^{22}(n))
+\frac{C_j^{12}(n)}{h_n^2}\sum_{k}\de_kC_k^{12}(n)
=-\frac{C_j^{12}(n)}{h_n}.\nonumber
\eq
Finally, using (6.18),
\bq
\pdj\ln\frac{\tau_{n+1}}{\tau_n}=-\frac{C_j^{12}(n)}{h_n}
=\pdj\ln h_n.
\eq
Summing over $n$ from $0$ to $N-1,$ we conclude that
$\tau_N$ is a constant multiple of $D_N,$ where the constant is
independent of $\{\de_j\}_{j=1}^{2g+2}.$ Since the $\tau$ - function
is defined up to such a constant, we can assume that the constant
is unity;
\bq
\tau_N(\de_1,...,\de_{2g+2})=D_N[w_{+}].
\eq

{\bf Remark 6.1.} It is worth mentioning that equations (\ref{freud1})
follow also (by putting $z = \delta_{j})$ from the equation
\begin{equation}\label{freud2}
A(z, n+1)U_n(z) - U_n(z)A(z, n) = \frac{\partial U_n(z)}{\partial z},
\end{equation}
which, in turn, is the compatibility condition of the basic
Fuchsian equation (\ref{fuchs}) and the difference equation
(\ref{difference}). This is the matrix form of the
so-called Freud equation which in principal can be written for any
semi-classical polynomials - see \cite{gan} and \cite{mag}
(and also \cite{FIK}).
In the physical language this is the ``discrete string equation''
corresponding to the weight $w_{+}(t)$. More precisely,
equation (\ref{freud2}) is the (discrete) Lax representation
of the Freud equation which manifests its integrability
from the algebraic point of view:
linear equations (\ref{fuchs}) and (\ref{difference}) form
a Lax pair for the Freud equation (cf. \cite{FIK}, \cite{FIK2}).
\vskip .5in
\setcounter{equation}{0}
{\section {Non-linear difference equations.}}

 As explained in Remark 6.1, the matrix equation (\ref{freud2})
should lead to the nonlinear difference equations for the recurrence coefficients,
following the {\it genre} of the Freud equations for the  Akhiezer polynomials.
To this end, we rewrite (6.24) elementwise,
by first specializing $\de_j$ to $\al_j$ and second to $\bt_j$. This will
produce six difference equations, relating polynomial evaluations at the
branch points and the recurrence coefficients.  For later convenience
we introduce four quantities
\bq
\mbox{} \ra_n&:=&\frac{1}{2h_{n-1}}P_{n}(\al_j)Q_{n-1}(\al_j),\nonumber\\
\mbox{} \rb_n&:=&\frac{1}{2h_{n-1}}P_{n}(\bt_j)Q_{n-1}(\bt_j),\nonumber\\
\mbox{} \Ra_n&:=&\frac{1}{2h_n}P_n(\al_j)Q_n(\al_j),\nonumber\\
\mbox{} \Rb_n&:=&\frac{1}{2h_n}P_n(\bt_j)Q_n(\bt_j).\nonumber
\eq
 Thus by specializing to $\al_j,$ $C_j(n)$ becomes,
\bq
\left(\matrix{\ra_n&-h_n\Ra_n\cr
              \Ra_{n-1}/h_{n-1}&-\ra_n-1/2\cr}\right),\nonumber
\eq
where we have taken onto account that the trace of the above is $-1/2.$ In
component form (6.24) is equivalent to,
\bq
\mbox{} \ra_{n+1}+\ra_n+\frac{1}{2}&=&\Ra_n(\al_j-b_{n+1})\\
\mbox{} a_{n+1}\Ra_{n+1}-a_n\Ra_{n-1}&=&(b_{n+1}-\al_j)
\left(\Ra_n(b_{n+1}-\al_j)+2\ra_n+\frac{1}{2}\right).
\eq
Note that out of the four possible equations, the $21$ element is a
tautology and the $11$ and $22$ elements are equivalent. Similarly,
specializing to $\bt_j,$
$C_j(n)$ becomes
\bq
\left(\matrix{\rb_n+\frac{1}{2}&-h_n\Rb_n\cr
              \Rb_{n-1}/h_{n-1}&-\rb_n\cr}\right),\nonumber
\eq
where the trace of the above is $1/2.$ In component form, (6.24) becomes,
\bq
\mbox{} \rb_{n+1}+\rb_{n}+\frac{1}{2}&=&\Rb_n(\bt_j-b_{n+1})\\
\mbox{} a_{n+1}\Rb_{n+1}-a_{n}\Rb_{n-1}&=&(\bt_{j}-b_{n+1})
\left(\Rb_n(\bt_j-b_{n+1})-2\rb_n-\frac{1}{2}\right).
\eq
In addition to these we have
\bq
\mbox{} a_n\Ra_n\Ra_{n-1}&=&\ra_n\left(\ra_n+\frac{1}{2}\right)\\
\mbox{} a_n\Rb_n\Rb_{n-1}&=&\rb_n\left(\frac{1}{2}+\rb_n\right),
\eq
since $\det C_j(n)=0.$ The equations (7.1) - (7.6)
are the difference equations mentioned above. We should be able to eliminate,
$\ra_n,\:\rb_n,\:\Ra_n$ and $\Rb_n$ from these to obtain non-linear
difference equations involving only $a_n$ and $b_n.$  These equations,
are also discussed in \cite{mag}.
\vskip .5cm
\setcounter{equation}{0}
{\section {\bf {The $\sigma_1$ Riemann-Hilbert Problem.}}}

In this section we shall solve the Riemann-Hilbert
problem {\bf{RH1}} - {\bf{RH4}} for the
Akhiezer polynomials in terms of the $\Theta$ -
functions.  To this end we will need a further
transformation of the Riemann-Hilbert problem satisfied
by $\Phi_n(z)$ to the so-called $\sigma_1$
problem, first  appeared in the theory of algebrogeometric solutions of
integrable PDEs (see \cite{its2},  \cite{bel}).

We notice that since the matrices $\left(\matrix{1&-1\cr
0&-1\cr}\right)$ and $\sigma_{1}$ have the same simple spectrum,
they must be similar. Indeed, we have
\bq
\left(\matrix{1&0\cr
              1&-1\cr}\right)\left(\matrix{1&-1\cr
                                           0&-1\cr}\right)
\left(\matrix{1&0\cr
              1&-1\cr}\right)=\left(\matrix{0&1\cr
                                            1&0\cr}\right)=\sigma_1.
\nonumber
\eq
Therefore, if we define
\bq \label{Psidef}
\mbox{} \Psi_n(z)&:=&\left(\matrix{\frac{1}{\sqrt{2\pi i}}&0\cr
                                   0&\sqrt{\frac{2}{\pi i}}\cr}\right)\Phi_n(z)
\left(\matrix{1&0\cr
     1&-1\cr}\right)\nonumber\\
\mbox{} &=&\left(\matrix{\frac{1}{\sqrt{2\pi i}}&0\cr
            0&\sqrt{\frac{2}{\pi i}}\cr}\right)Y_n(z)\left(\matrix{1&0\cr
                                                  0&1/w(z)\cr}\right)
\left(\matrix{\stp&0\cr
              1/\stp&-1/\stp\cr}\right),
\eq
then the jump matrix of the new function becomes $\sigma_{1}$. The left
diagonal constant matrix multiplier is introduced to normalize
the asymptotic behavior of the function $\Psi_n(z)$ at $z =\infty$:
$$
\left(\matrix{\frac{1}{\sqrt{2\pi i}}&0\cr
                                   0&\sqrt{\frac{2}{\pi i}}\cr}\right)\left(\matrix{z^n&0\cr
                                           0&z^{-n+1}\cr}\right)
\left(\matrix{\sqrt{2\pi i}&0\cr
                                   0&-\sqrt{\frac{\pi i}{2}}\cr}\right)
\left(\matrix{1&0\cr
              1&-1\cr}\right)
$$
$$
=\left(I+{\rm O}\left(\frac{1}{z}\right)\right)
\left(\matrix{z^{n}&0\cr
                     0&z^{-n+1}\cr}\right).
$$
Taking also into account that
$$
\left(\matrix{0&1\cr
              2&-1\cr}\right)\left(\matrix{1&0\cr
                                           1&-1\cr}\right)
= \left(\matrix{1&-1\cr
                1&1\cr}\right),
$$
we can reformulate the Riemann-Hilbert problem in terms of
$\Psi_n(z),$ as follows.
\bq
\mbox{} \Psi 1.&\quad&\Psi_n(z)\;{\rm is\;holomorphic\;for\;}
z\in\Co\setminus E.\nonumber\\
\mbox{} \Psi 2.&\quad&\Psi_{n-}(z)=\Psi_{n+}(z)\sigma_1,\quad z\in E.\nonumber
\\
\mbox{} \Psi 3.&\quad&\Psi_n(z)=\left(I+{\rm O}\left(\frac{1}{z}\right)\right)
z^{\left(\matrix{n&0\cr
                 0&-n+1}\right)},\quad z\rightarrow\infty.\nonumber\\
\mbox{} \Psi 4.&\quad&\Psi_n(z)=\hat{\Psi}_{n}^{(\beta_{j})}(z)
\left(\matrix{{\sqrt {z-\bt_j}}&0\cr
                              0&1\cr}\right)
\left(\matrix{1&-1\cr
              1&1\cr}\right)\nonumber\\
\mbox{} &\quad\quad=&\hat{\Psi}_{n}^{(\beta_{j})}(z)(z-\bt_j)^{\left(\matrix{1/2&0\cr
                                                            0&0\cr}\right)}
\left(\matrix{1&-1\cr
              1&1\cr}\right),\\
\mbox{} \Psi 5.&\quad&\Psi_n(z)=
\hat{\Psi}_{n}^{(\al_{j})}(z)
\left(\matrix{1/{\sqrt {z-\al_j}}&0\cr
                              0&1\cr}\right)
\left(\matrix{1&-1\cr
              1&1\cr}\right)\nonumber\\
\mbox{} &\quad\quad=&\hat{\Psi}_{n}^{(\al_{j})}(z)(z-\al_j)^{\left(\matrix{-1/2&0\cr
                                                            0&0\cr}\right)}
\left(\matrix{1&-1\cr
              1&1\cr}\right).\nonumber
 \eq
where $\hat{\Psi}_{n}^{(\al_{j})}(z)$ is holomorphic in the neighbourhood of
$z=\al_j$ and $\det \hat{\Psi}_{n}^{(\al_{j})}(\al_j)\neq 0,$i.e.,
\bq
\hat{\Psi}_{n}^{(\al_{j})}(z)=\sum_{k=0}^{\infty}\Psi_{nk}^{(\al_{j})}(z-\al_j)^k,
\quad\det\Psi_{n0}^{(\al_{j})}\neq 0.\nonumber
\eq
Similarly, $\hat{\Psi}_{n}^{(\bt_{j})}(z)$ is holomorphic in the neighbourhood
of $z=\bt_j$ and $\det \hat{\Psi}_{n}^{(\bt_{j})}(\bt_j)\neq 0,$i.e.,
\bq
\hat{\Psi}_{n}^{(\bt_{j})}(z)=\sum_{k=0}^{\infty}\Psi_{nk}^{(\bt_{j})}(z-\bt_j)^k,
\quad\det\Psi_{n0}^{(\bt_{j})}\neq 0.\nonumber
\eq
It is also worth noticing that the matrix products
$$\left(\matrix{{\sqrt {z-\bt_j}}&0\cr
                              0&1\cr}\right)\left(\matrix{1&-1\cr
1&1\cr}\right)$$ and $$\left(\matrix{1/{\sqrt {z-\al_j}}&0\cr
                              0&1\cr}\right)\left(\matrix{1&-1\cr
1&1\cr}\right)$$ have an exact $\sigma_{1}-$ jump matrix in the respective neighborhoods.

{\bf Remark 8.1.} From $\Psi 1 - \Psi 5$ it follows (independent of (\ref{Psidef})) that
\bq
\det\Psi_n(z)=\frac{i}{\pi w(z)},
\eq
{\bf Remark 8.2.} The function $\Psi_n(z),$ in terms of $P_n(z)$ and $Q_n(z),$ is
given as:
\bq\label{PsiPn}
\Psi_n(z)=\frac{1}{2\pi i}
\left(\matrix{\frac{i\pi w(z)P_n(z)-Q_n(z)}{w(z)}&
\frac{i\pi w(z)P_n(z)+Q_n(z)}{w(z)}\cr
2\frac{i\pi w(z)P_{n-1}(z)-Q_{n-1}(z)}{h_{n-1}w(z)}&
2\frac{i\pi w(z)P_{n-1}(z)+Q_{n-1}(z)}{h_{n-1}w(z)}\cr}\right),
\eq
and all the properties listed in $\Psi 1 - \Psi 5$ can be deduced from this
representation. It is worth emphasizing here that our approach does not
require this formula. Our logic is: The initial Riemann-Hilbert Problem
for $Y_n(z)$, quite generally posed, is transformed via (8.1) to the $\sigma_1$
problem which in turn leads to the equations (8.2) and (8.3)  by the
completely general principals of the Riemann-Hilbert problem.

 Let us now solve the $\sigma_1$ problem defined by $\Psi 1 - \Psi 5,$
however, without any reference to (8.4). The philosophy we adopt here is
similar to that in the asymptotic analysis of orthogonal polynomials via
the Riemann-Hilbert problem (cf. \cite{BI}, \cite{DKMVZ}): 
We simply ``forget'' the explicit formulas involving polynomials.

 Introduce the genus $g$ Riemann surface $\Rg$ defined by
\bq
y^2=(z-\bt_0)(z-\bt_{g+1})\prod_{j=1}^{g}(z-\al_j)(z-\bt_j),
\nonumber
\eq
and let ${\vP}_n(P),$ where $P=(z,y)\in \Rg$ be the vector Baker-Akhiezer
function determined by the conditions:
\bq
\mbox{} &{\bf {BA1}}.& \vP_n(P) {\;\rm is\;meromorphic\;on\;\;}\Rg\setminus
{\infty^{\pm}}\;{\rm with\;the\;pole\;divisor,\;}\nonumber\\
&\;\;\;&(\vP_n(P))=-\sum_{j=1}^{g}\al_j\nonumber\\
\mbox{} &{\bf {BA2}}.& {\rm The\;behaviour\;of\;}\vP_n(P){\;\rm at\;}
\infty^{\pm}
{\;\rm is\;specified\;by\;the\;equations,\;}\nonumber\\
\mbox{} &\;\;\;&\vP_n(P)=\left(\left(\matrix{1\cr
                                       0\cr}\right)+
{\rm O}\left(\frac{1}{z}\right)\right)z^n,\quad P\rightarrow\infty^{+},
\nonumber\\
\mbox{} &\;\;\;&\vP_n(P)=\left(\left(\matrix{0\cr
                                      1\cr}\right)+{\rm O}
\left(\frac{1}{z}\right)\right)z^{-n+1},\quad P\rightarrow\infty^{-},
\nonumber
\eq
in other words, $\infty^+$ is a pole of order $n$ and $\infty^{-}$ is a
zero of order $n-1.$ Here as usual, $\infty^{\pm}$ means
\bq
P\rightarrow\infty^{\pm}\iff z\rightarrow\infty,\;\;y\rightarrow\pm z^{g+1}.
\nonumber
\eq
Let ${\cal \pi}:\Rg\rightarrow \cpo$ be the projection,
\bq
{\cal \pi}(P)=z,\quad P=(z,y),\nonumber
\eq
and $*:\Rg\rightarrow\Rg^*$ be the involution,
\bq
P\rightarrow P^*=(z,-y)\;\;{\rm if\;\;}P=(z,y).\nonumber
\eq
The main observation (cf. \cite{its2},  \cite{bel}) is that the matrix function,
\bq\label{algPsidef}
\Psi_n(z):=\left(\vP_n(P),\vP_n(P^*)\right),
\eq
where
${\cal \pi}(P)=z,$ and $P\rightarrow\infty^+$ as $z\rightarrow\infty$,
solves the RH problem $\Psi 1 - \Psi 5.$

1. Indeed $\Psi 1$ is satisfied by construction since  (\ref{algPsidef}) defines $\Psi_n(z)$
uniquely as an analytic function on $\cpo\setminus E.$

2. If $z\rightarrow E$ from the ``+''-side (or from above the cut), then
\bq
\mbox{} P\rightarrow (z,y_{+}(z))=P_+\nonumber\\
\mbox{} P^*\rightarrow (z,-y_{+}(z))=(z,y_{-}(z))=P_{-}.\nonumber
\eq
If $z\rightarrow E$ from the ``-'' side, then
\bq
\mbox{} P\rightarrow (z,y_{-}(z))=P_{-}\nonumber\\
\mbox{} P^*\rightarrow(z,-y_{-}(z))=(z,y_{+}(z))=P_{+}.\nonumber
\eq
Hence,
\bq
\mbox{} \Psi_{n-}(z)&=&\left(\vP_n(P_{-}),\;\vP_n(P_{+})\right)\nonumber\\
\mbox{} \Psi_{n+}(z)&=&\left(\vP_n(P_{+}),\;\vP_n(P_{-})\right)\nonumber
\eq
and it follows that,
\bq
\Psi_{n-}(z)=\Psi_{n+}(z)\sigma_1,\quad z\in E,\nonumber
\eq
and therefore $\Psi 2$ is satisfied.

3. We have by construction, $z\rightarrow\infty$ implies
$P\rightarrow\infty^{+}$ and $P^*\rightarrow\infty^-$.

Therefore from ${\bf {BA2}},$
\bq
\Psi_n(z)=\left(I+{\rm O}\left(\frac{1}{z}\right)\right)
\left(\matrix{z^n&0\cr
              0&z^{-n+1}\cr}\right),
\eq
which shows that $\Psi 3$ is satisfied.

4. The function $\Psi(P)$ is analytic in the neighborhood of $P = \bt_{j}$
as a point of the Riemann surface $\Rg$. The local parameter
at the point $\bt_{j}$ is the square root of $z - \bt_{j}$.
Therefore, in the neighborhood of $P = \bt_{j}$ we have,
\bq
\mbox{} \vP_n(P)&=&\sum_{k=0}^{\infty}{\vec \psi}_{jk}(z-\bt_j)^{k/2},
\label{tylbt1}\\
\mbox{} \vP_n(P^*)&=&\sum_{k=0}^{\infty}(-1)^{k}{\vec \psi}_{jk}
(z-\bt_j)^{k/2},\label{tylbt2}
\eq
so that
$$
\Psi_{n}(z) =\left(\sum_{k=0}^{\infty}{\vec \psi}_{jk}(z-\bt_j)^{k/2},
 \sum_{k=0}^{\infty}(-1)^{k}{\vec \psi}_{jk}
(z-\bt_j)^{k/2}\right).
$$
This in turn implies that the function $\hat{\Psi}_{n}^{(\beta_{j})}(z)$
defined by  the equation $\Psi 4$ is {\it a holomorphic function of $z$}.
Indeed we have
$$
\hat{\Psi}_{n}^{(\beta_{j})}(z) \equiv \Psi_n(z)
\left(\matrix{1&-1\cr
              1&1\cr}\right)^{-1}
\left(\matrix{\frac{1}{\sqrt {z-\bt_j}}&0\cr
                              0&1\cr}\right)
$$
$$
=
\frac{1}{2}\left(\sum_{k=0}^{\infty}\left[{\vec \psi}_{jk}- (-1)^{k}
{\vec\psi}_{jk}\right](z-\bt_j)^{{k-1}/2},
 \sum_{k=0}^{\infty}\left[{\vec \psi}_{jk} + (-1)^{k}{\vec \psi}_{jk}\right]
(z-\bt_j)^{k/2}\right)
$$
$$
=\left(\sum_{l=0}^{\infty}{\vec \psi}_{j2l+1}(z-\bt_j)^{l},
 \sum_{l=0}^{\infty}{\vec \psi}_{j2l}
(z-\bt_j)^{l}\right).
$$

5. Since $P=\al_j$ is a simple pole of $\Psi(P)$, the Taylor
series (\ref{tylbt1}) and (\ref{tylbt2}) shoud be replaced by
the Laurent series,
\bq
\mbox{} \vP_n(P)&=&\sum_{k=-1}^{\infty}{\vec \phi}_{jk}(z-\al_j)^{k/2},
\nonumber\\
\mbox{} \vP_n(P^*)&=&\sum_{k=-1}^{\infty}(-1)^{k}{\vec \phi}_{jk}
(z-\al_j)^{k/2}.\nonumber
\eq
The rest of the arguments is literaly the same as in the $\bt$ -case,
and we have that the function $\hat{\Psi}_{n}^{(\al_{j})}(z)$
defined by  the equation $\Psi 5$ is holomorphic  at $z = \al_{j}$.

Our final observation is that already established properties imply
(8.3) (cf. our ``Riemann-Hilbert'' proof of (\ref{dety}) above)
and hence the inequalities,
$$
\det \hat{\Psi}_{n}^{(\alpha_{j})}(\alpha_{j}) \neq 0,
\quad \det \hat{\Psi}_{n}^{(\beta_{j})}(\bt_{j}) \neq 0.
$$

We now come to the $\Theta-$ formula for $\vP_n(P).$ First we assemble here
for this purpose some facts about the Riemann surface $\Rg$ realized 
as a two-sheet covering of the $z$ plane in the usual way and with the
first homology basis depicted in the figure below.
\begin{figure}[h]
\centering
\includegraphics[height=1.8in]{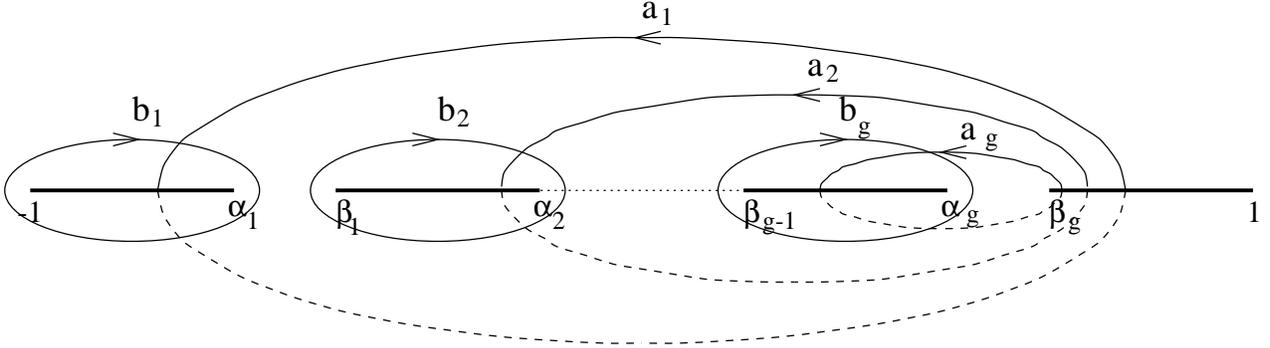}
\caption{The dash curves represent the parts of the cannonical loops lying
on the lower sheet. The lower (upper) sheet is fixed by the condition that
it contains the point $\infty^+$ ( $\infty^-$).}
\end{figure}
 Let  
$$
\{d\om_j\}_{j=1}^{g}, \quad \int_{a_j}d\om_k=\de_{jk},
$$ 
be a set of normalised Abelian differentials of the first kind.
As it is usual for a hyperelliptic curve, we shall chose the
differentials $d\om_j$ according to the equations,
\bq
\mbox{} d\om_j&=&\sum_{k=1}^{g}(A^{-1})_{jk}\frac{z^{g-k}}{y}dz,\nonumber\\
\mbox{} A_{jk}&=&\int_{a_k}\frac{z^{g-j}}{y}dz.\nonumber
\eq
The invertability of the matrix $A$ is a (relatively simple) classical result.
We refer the reader to the  monograph \cite{fkra} for the basic general facts
concerning the theory of functions on the Riemann surfaces (see also
chapter 1 of \cite{bel}). Let us also introduce the normalized Abelian differential 
of the third kind, having its only poles
at $\infty^{\pm},$
\bq
d\Om(P)=\frac{z^g+\lambda_{g-1}z^{g-1}+...+\lambda_0}{y}dz,\nonumber
\eq
with vanishing $a-$period;
\bq
\int_{a_j}d\Om=0,\quad j=1,...,g.\nonumber
\eq
The above $g$ conditions uniquely determine \cite{fkra} the coefficients,
$\{\lambda_j\}_{j=0}^{g-1}$. Put
$$
\Omega(P) = \int_{\beta_{g+1}}^{P}d\Om.
$$ 
One easily deduces,
\bq \label{Omasymp}
\Om(P)=\pm\left(\ln z - \ln C(E) + {\rm O}\left(\frac{1}{z}\right)\right),\quad
P\rightarrow\infty^{\pm},
\eq
where
\bq\label{diam}
C(E)=\exp\left(-\int_{\bt_{g+1}}^{\infty^{+}}
\left(\frac{z^g+\sum_{j=0}^{g-1}\lambda_jz^j}{y(z)}-\frac{1}{z}\right)
dz\right).
\eq
(We recall that $\beta_{j+1} = 1$.)
Finally, the Riemann $\Theta-$function of $g-$ complex variables
${\vec s}\in\Co^g,$ is defined with the aid of the period matrix
\bq
B_{jk}:=\int_{b_k}d\om_j,\nonumber
\eq
as follows:
\bq
\Th({\vec s})\equiv \Th({\vec s};B):=\sum_{{\vec t}\in \cz^g}\exp\left(\:i\pi({\vec t},B{\vec
t}\:) +2\pi i({\vec t},{\vec s}\:)\:\right).\nonumber
 \eq
Here are the fundamental periodic  property of the $\Theta-$ function:
\begin{equation}\label{thetaperiodic}
\Th({\vec s} + \vec n + B\vec m)
= e^{-\pi i(B\vec m, \vec m) - 2\pi i (\vec s, \vec m)}
\Th({\vec s}),
\end{equation}
and the obvious symmetry relation:
$$
\Th({-\vec s}) = \Th({\vec s}).
$$
\vskip.2in
Observe now that ${\bf {BA1 - BA2}}$
imply the following properties on the components
of $\vP_n(P).$ 
\bq
\mbox{} &\Psi_{n1}(P)&\;{\;\rm is\;meromorphic\;on\;}\Rg
\setminus\{\infty^{+},\infty^{-}\}\nonumber\\
\mbox{} &(\Psi_{n1}(P))&=-\sum_{j=1}^{g}\al_j\nonumber\\
\mbox{} &\Psi_{n1}(P)&=z^n+{\rm O}(z^{n-1}),\quad P\rightarrow\infty^{+}
\nonumber\\
\mbox{} &\Psi_{n1}(P)&={\rm O}(z^{-n}),\quad P\rightarrow\infty^{-}.\label{aspsi1}
\eq
Similary for $\Psi_{n2}(P),$
\bq
\mbox{} \Psi_{n2}(P)&=&z^{-n+1}+{\rm O}(z^{-n}),\quad
O\rightarrow\infty^{-},\nonumber\\
\mbox{} \Psi_{n2}(P)&=&{\rm O}(z^{n-1}),\quad P\rightarrow\infty^{+}.\label{aspsi2}
\eq
By standard technique of the algebrogeometric method (
see e.g. \cite{bel}), we get,
\bq
\mbox{} \Psi_{n1}(P)&=&{\rm e}^{n\Om(P)}
\frac{\Th\left(\int_{\bt_{g+1}}^{P}d\vo+n\vb-\vd\right)}
{\Th\left(\int_{\bt_{g+1}}^{P}d\vo-\vd\right)}
\frac{\Th\left(\int_{\bt_{g+1}}^{\infty^+}d\vo-\vd\right)}
{\Th\left(\int_{\bt_{g+1}}^{\infty^+}\dvo+n\vb-\vd\right)}C^{n}(E),\nonumber\\
\mbox{} \Psi_{n2}(P)&=&{\rm e}^{(n-1)\Om(P)}
\frac{\Th\left(\int_{\bt_{g+1}}^{P}d\vo+(n-1)\vb-\vd\right)}
{\Th\left(\int_{\bt_{g+1}}^{P}d\vo-\vd\right)}
\frac{\Th\left(\int_{\bt_{g+1}}^{\infty^+}d\vo+\vd\right)}
{\Th\left(\int_{\bt_{g+1}}^{\infty^+}\dvo-(n-1)\vb+\vd\right)}C^{(1-n)}(E),\nonumber
\eq
where
\bq
\mbox{} L_j&=&\frac{1}{2\pi i}\int_{bj}d\Om\nonumber\\
\mbox{} D_j&=&\sum_{k=1}^{g}\int_{\bt_{g+1}}^{\al_k}d\om_j+C_j\nonumber\\
\mbox{} &=&2\sum_{k=1}^{g}\int_{\bt_{g+1}}^{\al_k}d\om_j,\nonumber
\eq
and $C_{j}$ form the vector of the Riemann constants (see again \cite{fkra} and
\cite{bel}). Indeed, by the Riemann theorem (see e.g. \cite{fkra}), 
the first $\Theta-$functions in the denominators has zeros exactly 
at the points $\alpha_{j}$; the front exponential factors provide 
the needed asymptotic behavior at $\infty ^{\pm}$; the first 
$\Theta-$ functions in the numerators, by virture of the periodicity 
property (\ref{thetaperiodic}),
ensure the single-valuedness; and, finally, the $P$-independent
$\Theta$-factors
together with the back exponential factors provide the needed normalizations
ad  $\infty ^{\pm}$ (cf. (\ref{aspsi1}) and (\ref{aspsi2})). We also
assume that we choose the same path between $\beta_{g+1}$ and $P$  for
all the integrals involved{\footnote{Alternatively, one can choose for each
integral its own path. In this case though the paths must not intersect
the basic cycles.}}.

The formulae above can be simplified. To this end we observe that
\begin{equation}\label{halfperiod}
\int_{\beta_{g+1}}^{\alpha_{k}} d\omega_{j} =
\frac{1}{2}\delta_{jk} + \frac{1}{2}\sum_{l=1}^{k}B_{jl},
\end{equation}
where the path of integration from  $\beta_{g+1}$ to 
$\alpha_{k}$ lies on the upper plane of the upper sheet. 
Therefore, moduli the lattice periods, 
$$
D_{j} = 1 + \sum_{k=1}^{g}B_{jk}(g-k+1).
$$
In other words, the vector $\vd$ belongs to the latice 
$\cz^g + B\cz^g$ and hence (property (\ref{thetaperiodic}) again) 
can be droped from the above formulae for $\vP_n(P)$. This yields the following
simplified  $\Theta-$ representation for $\vP_n(P).$
\bq
\mbox{} \Psi_{n1}(P)&=&{\rm e}^{n\Om(P)}
\frac{\Th\left(\int_{\bt_{g+1}}^{P}d\vo+n\vb\right)}
{\Th\left(\int_{\bt_{g+1}}^{P}d\vo\right)}
\frac{\Th\left(\int_{\bt_{g+1}}^{\infty^+}d\vo\right)}
{\Th\left(\int_{\bt_{g+1}}^{\infty^+}\dvo+n\vb\right)}C^{n}(E),\label{Psitheta1}\\
\mbox{} \Psi_{n2}(P)&=&{\rm e}^{(n-1)\Om(P)}
\frac{\Th\left(\int_{\bt_{g+1}}^{P}d\vo+(n-1)\vb\right)}
{\Th\left(\int_{\bt_{g+1}}^{P}d\vo\right)}
\frac{\Th\left(\int_{\bt_{g+1}}^{\infty^+}d\vo\right)}
{\Th\left(\int_{\bt_{g+1}}^{\infty^+}\dvo-(n-1)\vb\right)}C^{(1-n)}(E),\nonumber
\eq

We conclude the $\Theta-$ function solution of the Akhiezer Riemann-Hilbert
problem by noticing the following equation for the vector $\vb$ of the
$b$ - periods of the integral $\Omega(P)$.
 
$$
{\vb} = \mbox{res}_{P = \infty^+}(\vec{\omega}d\Om(P))
\, + \,\mbox{res}_{P = \infty^-}(\vec{\omega}d\Om(P))
$$
\begin{equation}\label{Bom}
= -\int_{\bt_{g+1}}^{\infty^+}\dvo+
\int_{\bt_{g+1}}^{\infty^-}\dvo\, \, 
=\, \, -2\int_{\bt_{g+1}}^{\infty^+}\dvo,
\end{equation}
The equation  is just the classical Riemann bilinear 
identity (see e.g. \cite{fkra} or \cite{bel})
applied to the pair  of the Abelian integrals 
$\vec{\omega}(P)$ and $\Omega(P)$.

\vskip .2in

{\bf Remark 8.3} Using equation (\ref{halfperiod}), one can
check directly, with the help of the periodic condition
(\ref{thetaperiodic}), that the theta function,
$$
\Th\left(\int_{\bt_{g+1}}^{P}d\vo\right)
$$
has the points $\alpha_{j}$ as its zeros. 

\vskip .2in

{\bf Remark 8.4} The reader should not be misled by the formal 
possibility to diagonalize simultaneously  the jump matrices of the 
Riemann-Hilbert problem  $\Psi 1 - \Psi 5$
(which all are equal to $\sigma_{1}$) and by apparently following
from this conclusion that the problem can be reduced 
to the scalar one on the complex plane and hence solved
without any use of the $\Theta-$ functions. The obstractions come
from the end points $\alpha_{j}$, $\beta_{j}$ and from the
point at infinity, where the function $\Psi_{n}(z)$
must have the singularities specified by equations
$\Psi 5$, $\Psi 4$ and  $\Psi 3$, respectively. These singularities
can be alternatively discribed as the addition {\it jump
conditions} posed on the small circles around the end points
and on the big circle around the infinity. The relevant jump
matrices are
$$
\left(\matrix{{1/\sqrt {z-\alpha_j}}&0\cr
                              0&1\cr}\right)
\left(\matrix{1&-1\cr
              1&1\cr}\right), \quad
\left(\matrix{{\sqrt {z-\bt_j}}&0\cr
                              0&1\cr}\right)
\left(\matrix{1&-1\cr
              1&1\cr}\right),
$$
and
$$
z^{\left(\matrix{n&0\cr
                 0&-n+1}\right)},
$$
respectively. Posed in this form, the $\sigma_{1}$
Riemann-Hilbert problem becomes the regular one - no 
singularities different from the jumps are prescribed.
At the same, the  additional jump matrices depend on $z$ and
the whole new set of jump matrices can not be simultaneously 
diagonalized. The only way to circumvent this obstacles,
and not to use the $\Theta-$ functions, is the equation
(8.4) which indeed gives an explicit representation
of the solution of the $\sigma_{1}$ Riemann-Hilbert problem
in terms of the elementary functions and their contour
integrals. The $\Theta-$ function representation
(\ref{Psitheta1}) for
the solution $\Psi_{n}(z)$ obtained in this chapter has an important
advantage comparing to (8.4). It reveals the nature of
the dependence of  $\Psi_{n}(z)$, and hence of the Akhiezer
polynomials themselves (see (\ref{Pntheta}) below), on the
number $n$, as $n$ varies over the whole range $1 \leq n \leq \infty$
(see \cite{Chen} for more on the use of the $\Theta$ - representations
in the analysis of the  Akhiezer polynomials). Simultaneously,
the comparison of equations (8.4) and (\ref{Psitheta1}) might,
perhaps, be used to derive some new nontrivial identities 
for the hyperelliptic $\Theta-$ functions.
\vskip .2in

{\bf Remark 8.5} Up to a trivial diagonal gauge transformation,
the matrix function $\Psi_{n}(z)$ satisfies the same Fuchsian
equation (\ref{fuchs}) that is satisfied by the function 
$\Phi_{n}(z)$. Note that the corresponding monodromy group
is very simple; indeed, it has just one generator - the 
matrix $\sigma_{1}$. Once again, the reader might be
wondering about the appearance of the highly nontrivial
theta-functional formulae in the describtion of the
function $\Psi_{n}(z)$ which gives the
solution of the corresponding inverse monodromy problem.
Similar to the previous remark, the explanation 
comes from the fact that the solution
$\Psi_{n}(z)$, in addition to the given monodromy group,
must exhibit the local behavior at the singular points
indicated by the conditions  $\Psi 3 - \Psi 5$.
This situation is typical in the theory of the finite-gap
solutions of integrable PDEs
{\footnote{Another example of an apparently simple
but nontrivialy solved invesre monodromy problem can
be also found in the theory
of integrable PDEs. It is provided
by the multi-soliton Baker-Akhiezer function whose monodromy
group is just trivial. Of course, 
the formulae in this case are simplier than the
finite-gap ones - they do not
contain the $\Theta-$ functions. At the same time, the answer
is still rather complicated; in fact, 
it involves degenerated $\Theta-$ functions corresponding
to the singular curves of genus zero.}}
(see e.g. \cite{jmu} and \cite{bel}).

\setcounter{equation}{0}
{\section {\bf {A list of the $\Theta$ - formulae.}}

In this section, we give formulae expressing the polynomial $P_{n}(z)$, 
recurrence coefficients $a_n,$ $b_n,$ the square of the weighted $L^2$ norm $h_n$ and
the Hankel determinant in terms of the
$\Theta-$ functions. The expressions will be derived as simple corollaries
of the equations (\ref{algPsidef}) and (\ref{Psitheta1}) representing the
solution $\Psi_{n}(z)$ of the Riemann-Hilbert problem  $\Psi 1 - \Psi 5$
in terms of the $\Theta$ - functions.

From (\ref{Psidef}) it follows that (see also (\ref{PsiPn}))
$$
P_{n}(z) = (Y_{n}(z))_{11} = (\Psi_{n}(z))_{11} + (\Psi_{n}(z))_{12}.
$$
This together with  (\ref{algPsidef}) and (\ref{Psitheta1}) leads to
the following $\Theta$ - representation of the Akhiezer polynomials,

$$
P_{n}(z) = 
\frac{\Th\left(n\vb + \int_{\bt_{g+1}}^{z}d\vo\right){\rm e}^{n\Om(z)} + 
\Th\left(n\vb - \int_{\bt_{g+1}}^{z}d\vo\right){\rm e}^{-n\Om(z)}}
{\Th\left(\int_{\bt_{g+1}}^{z}d\vo\right)}
$$

\begin{equation}\label{Pntheta}
\times\frac{\Th\left(\int_{\bt_{g+1}}^{\infty^+}d\vo\right)}
{\Th\left(\int_{\bt_{g+1}}^{\infty^+}\dvo+n\vb\right)}C^{n}(E),
\end{equation}
\vskip.2in
\noindent
where all the hyperelliptic integrals are taken in the upper sheet
of the curve $\Rg$ (and along the same path). 

{\bf Remark 9.1} It is a simple but an instructive exercise to check directly,
using equation (\ref{halfperiod}), the similar equation for
the integral $\Omega(P)$, i.e.
$$
\Omega(\alpha_{k}) =
\pi i + \pi i \sum_{j=1}^{k}L_{j},
$$
and, once again, the periodicity property of the $\Theta$-function, that
the right side of (\ref{Pntheta}) is indeed a {\it polynomial}.

To  evaluate the quantities $a_n,$ $b_n,$ and $h_n$
we shall use the relation 
\begin{equation}\label{psi1m1}
\psi_{1} =\left(\matrix{\frac{1}{\sqrt{2\pi i}}&0\cr
              0&\sqrt{\frac{2}{\pi i}}\cr}\right) m_{1}
\left(\matrix{\sqrt{2\pi i}&0\cr
              0&\sqrt{\frac{\pi i}{2}}\cr}\right) 
- \left(\matrix{0&0\cr
              0&\kappa\cr}\right), \quad n > 1, 
\end{equation}
between the first matrix coefficients, $\psi_{1}$ and  $m_{1}$,
of the  Laurent series
$$
\Psi_n(z)=\left(I+\sum_{k=1}^{\infty}\frac{\psi_k(n)}{z^k}\right)
\left(\matrix{z^{n}&0\cr
              0&z^{-n +1}\cr}\right),
\quad |z| > 1,
$$
and
$$
Y_n(z)=\left(I+\sum_{k=1}^{\infty}\frac{m_k(n)}{z^k}\right)z^{n\st},
\quad |z| > 1,
$$
respectively. In (\ref{psi1m1}), the  parameter $\kappa$ is defined
via the expansion,
$$
w(z) = \frac{i}{\pi z}\left( 1 + \frac{\kappa}{z} + ...\right),
$$ 
and the matrix,
$$
\left(\matrix{0&0\cr
              -1&0\cr}\right),
$$
should be added to the r.h.s if $n =1$. Combaining equation (\ref{psi1m1})
with the formula (\ref{hab_Y1}) we obtain that
$$
h_{n} = 2\left(\psi_{1}(n)\right)_{12}.
$$
On the other hand, let us introduce the coefficient matrix
$c_{jk}$, $j,k = 1, 2 $ by the relations (cf. (\ref{aspsi1}) and (\ref{aspsi2})),
\bq
\mbox{} &\Psi_{n1}(P)&=z^n+ c_{11}z^{n-1} + {\rm O}(z^{n-2}),
\quad P\rightarrow\infty^{+}
\label{c11}\\
\mbox{} &\Psi_{n1}(P)&=c_{12}z^{-n} + {\rm O}(z^{-n-1}),
\quad P\rightarrow\infty^{-},\label{c12}
\eq
and 
\bq
\mbox{} \Psi_{n2}(P)&=&z^{-n+1}+ c_{22}z^{-n} + {\rm O}(z^{-n - 1}),\quad
O\rightarrow\infty^{-},\label{c22}\\
\mbox{} \Psi_{n2}(P)&=&c_{21}z^{n-1} + {\rm O}(z^{n-2}),\quad P\rightarrow\infty^{+}.
\eq
Then, it is obvious that 
\begin{equation}\label{psi1c}
\left(\psi_{1}(n)\right)_{jk}  = c_{jk},
\end{equation}
and, in particular, we arrive to the equation
\begin{equation}\label{hnc12}
h_{n} = 2c_{12}.
\end{equation}
The coefficient $c_{12}$, in its turn, can be immediately evaluated from
the  $\Theta -$ formula (\ref{Psitheta1}) by letting
$P \to \infty^{-}$. In fact, we have
\bq
\mbox{} c_{12}&=&C^{2n}(E)
\frac{\Theta\left(\int_{\bt_{g+1}}^{\infty^+}\dvo-n{\vb}\right)}
{\Theta\left(\int_{\bt_{g+1}}^{\infty^+}\dvo+n{\vb}\right)}.\label{hntheta}
\eq
Taking into account the Riemann bilinear relation (\ref{Bom}) we can 
present the formula for $h_{n}$ in the following final form,
\bq
\mbox{} h_{n}&=&2C^{2n}(E)
\frac{\Theta\left(\left(n + \frac{1}{2}\right){\vb}\right)}
{\Theta\left(\left(n - \frac{1}{2}\right){\vb}\right)},\quad
n=1,2,...,\label{hntheta1}\\
\mbox{} h_0&:=&1.\nonumber
\eq
An important direct consequence of this equation is the explicit
$\Theta-$ functional representation for determinant of 
the $(n+1)\times(n+1)$ Hankel matrix:
\bq
D_{n+1}[w_+]=\prod_{j=0}^{n}h_j=2^{n}(C(E))^{n(n+1)}
\frac{\Theta\left(\left(n + \frac{1}{2}\right){\vb}\right)}
{\Theta\left(\frac{1}{2}{\vb}\right)}\nonumber
\eq 
\begin{equation}\label{hankeltheta}
=2^{n}(C(E))^{n(n+1)}
\frac{\Theta\left((2n+1)\int_{\bt_{g+1}}^{\infty^+}\dvo\right)}
{\Theta\left(\int_{\bt_{g+1}}^{\infty^+}\dvo\right)}.
\end{equation}

A similar use of the remaining equations in (\ref{psi1m1}), (\ref{psi1c}) and
the formulae (\ref{hab_Y2}),  (\ref{hab_Y3}) leads
at once to the $\Theta-$ representations of the recurrence
coefficients $a_{n}$ and $b_{n}$:
\bq
a_n =
\left\{ \begin{array}{ll}
2C^2(E)\frac{\Theta\left(\frac{3}{2}\vb\right)}
{\Theta\left(\frac{1}{2}\vb\right)} &
 \textrm{if $n=1$}\\\\
C^2(E)\frac{\Theta\left((n +\frac{1}{2})\vb\right)
\Theta\left((n -\frac{3}{2})\vb\right)}
{\Theta^{2}\left((n-\frac{1}{2})\vb\right)} &
\textrm{if $n>1$}\end{array}\right.,\label{antheta}
\eq
and
\bq
b_n &=& \frac{1}{2}\sum_{j=1}^g(\bt_j-\al_j) \nonumber \\
&+&\sum_{j=1}^g(A^{-1})_{j1}
\left[{\scriptstyle\frac{\Theta_j^{\pr}\left((n-\frac{1}{2})\vb\right)}
{\Theta\left((n-\frac{1}{2})\vb\right)}-
\frac{\Theta_j^{\pr}\left((n-\frac{3}{2})\vb\right)}
{\Theta\left((n-\frac{3}{2})\vb\right)} -
2\frac{\Theta_j^{\pr}\left(\frac{1}{2}\vb\right)}
{\Theta\left(\frac{1}{2}\vb\right)}}\right].
\label{bntheta}
\eq
Here,
\bq
\Theta_{j}^{\pr}\left({\vec s}\right):=
\frac{\partial\Theta({\vec s})}{\partial s_j}\nonumber.
\eq
Equations (\ref{Pntheta}), (\ref{hntheta1}), (\ref{hankeltheta}),
(\ref{antheta}) and  (\ref{bntheta}) were previously obtained in
\cite{Chen} by a direct analysis of Akhiezer's 
function defined as the sum $\frac{i\pi w(z)P_n(z)-Q_n(z)}{w(z)}$
(cf. (\ref{PsiPn})).  In \cite{Chen} it was also shown that
the above formulae allow to identify the quantity $C(E)$ as
the transfinite diameter of the set $E$. We remind that in
our approach, $C(E)$ appears as a first nontrivial
coefficient in the asymptotic expansion of the 
Abelian integral $\Omega(P)$, see (\ref{Omasymp}) and (\ref{diam}). 
Finally, we should note that equations (\ref{Psitheta1}),  
(\ref{antheta}) and  (\ref{bntheta}), as the equations describing
the eigenfunctions and the coefficients of a finite-gap discrete
Schr\"dinger operator, have already
been known ( see e.g. \cite{krich}) in the theory of the 
periodic Toda lattice.

\end{document}